\documentclass{article}
\usepackage{graphicx}
\usepackage{tikz}
\usepackage{geometry}
\usepackage{setspace}
\usepackage{indentfirst}
\usepackage{xcolor}
\usepackage{hyperref}
\usepackage{amsmath}
\usepackage{amsfonts}
\usepackage{float}
\usepackage{amssymb,amsthm,graphicx}
\usepackage{tcolorbox} 

\hypersetup{
    colorlinks=true,
    linkcolor=blue,
    citecolor=blue,
    filecolor=blue,
    urlcolor=blue,
}

\title{\textbf{Topology-Preserving Scaling in Data Augmentation} \\[0.5em]}
\author{
\textbf{Vu Anh Le\textsuperscript{1*}, Mehmet Dik\textsuperscript{1,2**}}  \\
\textsuperscript{1} Department of Mathematics and Computer Science, Beloit College \\
\textsuperscript{2} Department of Mathematics, Computer Science \& Physics, Rockford University \\ [1em]
\textit{Contact} \\
* \href{mailto:csplevuanh@gmail.com}{csplevuanh@gmail.com} \\
** \href{mailto:mdik@rockford.edu}{mdik@rockford.edu}
}
\date{\today}

\begin{document}

\maketitle

\begin{abstract}
We propose an algorithmic framework for dataset normalization in data augmentation pipelines that preserves topological stability under non-uniform scaling transformations. Given a finite metric space \( X \subset \mathbb{R}^n \) with Euclidean distance \( d_X \), we consider scaling transformations defined by scaling factors \( s_1, s_2, \ldots, s_n > 0 \). Specifically, we define a scaling function \( S \) that maps each point \( x = (x_1, x_2, \ldots, x_n) \in X \) to
\[
S(x) = (s_1 x_1, s_2 x_2, \ldots, s_n x_n).
\]
Our main result establishes that the bottleneck distance \( d_B(D, D_S) \) between the persistence diagrams \( D \) of \( X \) and \( D_S \) of \( S(X) \) satisfies:
\[
d_B(D, D_S) \leq (s_{\max} - s_{\min}) \cdot \operatorname{diam}(X),
\]
where \( s_{\min} = \min_{1 \leq i \leq n} s_i \), \( s_{\max} = \max_{1 \leq i \leq n} s_i \), and \( \operatorname{diam}(X) \) is the diameter of \( X \). Based on this theoretical guarantee, we formulate an optimization problem to minimize the scaling variability \( \Delta_s = s_{\max} - s_{\min} \) under the constraint \( d_B(D, D_S) \leq \epsilon \), where \( \epsilon > 0 \) is a user-defined tolerance.

We develop an algorithmic solution to this problem, ensuring that data augmentation via scaling transformations preserves essential topological features. We further extend our analysis to higher-dimensional homological features, alternative metrics such as the Wasserstein distance, and iterative or probabilistic scaling scenarios. Our contributions provide a rigorous mathematical framework for dataset normalization in data augmentation pipelines, ensuring that essential topological characteristics are maintained despite scaling transformations.
\end{abstract}

\tableofcontents

\vspace{0.3in}

\section{Introduction}

Data augmentation is a popular technique in machine learning for enhancing model generalization by artificially increasing the diversity of training data. Common augmentation methods include geometric transformations such as rotations, translations, and scaling \cite{shorten2019survey}. In particular, scaling transformations are widely used due to their simplicity and effectiveness \cite{zhang2017mixup}. However, non-uniform scaling, where each coordinate axis is scaled by a distinct factor, can introduce anisotropic distortions that significantly alter the intrinsic geometry and topology of datasets \cite{chen2020group}.

Topological Data Analysis (TDA) provides an approach to capture the intrinsic shape of data in a way that is robust to noise and deformation \cite{edelsbrunner2010computational}. Critical to TDA is the concept of \textit{persistent homology}, which summarizes topological features of data across multiple scales using \textit{persistence diagrams} \( D \). A key property of persistence diagrams is their stability under perturbations of the input data, as quantified by the bottleneck distance \( d_B \) \cite{cohen2007stability}.

In prior work \cite{le2024stability}, we have investigated the effects of non-uniform scaling transformations defined by:
\[
S(x) = (s_1 x_1, s_2 x_2, \ldots, s_n x_n),
\]
where \( s_i > 0 \) for all \( i \). Our primary goal was to establish explicit bounds on the bottleneck distance \( d_B(D, D_S) \) between the persistence diagrams before and after scaling. Specifically, we showed that:
\[
d_B(D, D_S) \leq \delta = \frac{1}{2} \Delta_s \cdot \operatorname{diam}(X),
\]
where \( \Delta_s = s_{\max} - s_{\min} \). This inequality provides a direct relationship between the scaling variability \( \Delta_s \) and the topological perturbation measured by \( d_B(D, D_S) \).

Based on this theoretical guarantee, we formulate an optimization problem to minimize \( \Delta_s \) under the constraint \( d_B(D, D_S) \leq \epsilon \), where \( \epsilon > 0 \) is a user-defined tolerance. The solution to this problem yields scaling factors that minimize anisotropic distortions while preserving the topological features of the dataset. We further extend our analysis to consider higher homology dimensions \cite{otter2017roadmap}, alternative distance metrics such as the Wasserstein distance \cite{turner2014frechet}, and scenarios involving iterative or probabilistic scaling \cite{le2024stability}.

\vspace{0.3in}

\section{Preliminaries}

\subsection{Metric Spaces and Scaling Transformations}

Let \( X = \{ x_1, x_2, \ldots, x_N \} \subset \mathbb{R}^n \) be a finite metric space with the Euclidean distance \( d_X: X \times X \to \mathbb{R} \), defined by:
\[
d_X(p, q) = \| p - q \|_2 = \left( \sum_{i=1}^n (p_i - q_i)^2 \right)^{1/2}.
\]

Consider a scaling transformation \( S: \mathbb{R}^n \to \mathbb{R}^n \) given by:
\[
S(x) = (s_1 x_1, s_2 x_2, \ldots, s_n x_n),
\]
where \( s_i > 0 \) for \( 1 \leq i \leq n \). The scaled dataset is \( S(X) = \{ S(x) \mid x \in X \} \), and the scaled distance \( d_S \) between points \( p, q \in X \) is:
\[
d_S(p, q) = \| S(p) - S(q) \|_2 = \left( \sum_{i=1}^n s_i^2 (p_i - q_i)^2 \right)^{1/2}.
\]

\vspace{0.15in}

\subsection{Persistence Diagrams and Bottleneck Distance}

A \textit{filtration} \( \{ K_\epsilon \}_{\epsilon \geq 0} \) is a nested sequence of simplicial complexes built on \( X \), such that \( K_{\epsilon} \subseteq K_{\epsilon'} \) whenever \( \epsilon \leq \epsilon' \). Common filtrations include the Vietoris--Rips and \v{C}ech complexes.

The \textit{persistent homology} of \( X \) captures the birth and death times of topological features (e.g., connected components, loops, voids) as the scale parameter \( \epsilon \) varies. The collection of these features is summarized in the \textit{persistence diagram} \( D \), which is a multiset of points \( (b, d) \in \mathbb{R}^2 \), where \( b \) is the birth time and \( d \) is the death time of a feature.

The \textit{bottleneck distance} \( d_B(D, D') \) between two persistence diagrams \( D \) and \( D' \) is defined as:
\[
d_B(D, D') = \inf_{\gamma} \sup_{x \in D} \| x - \gamma(x) \|_{\infty},
\]
where \( \gamma: D \to D' \) is a bijection (allowing for matching points to the diagonal \( b = d \)), and \( \| \cdot \|_{\infty} \) denotes the \( L^\infty \)-norm.

\vspace{0.15in}

\subsection{Stability of Persistence Diagrams}

The stability theorem \cite{cohen2007stability} states that small perturbations in the input data lead to small changes in the persistence diagrams. Specifically, for two functions \( f, g: X \to \mathbb{R} \), the bottleneck distance between their persistence diagrams satisfies:
\[
d_B(D_f, D_g) \leq \| f - g \|_{\infty}.
\]

When considering metric spaces, if \( d_X \) and \( d_{X'} \) are distance functions on \( X \) satisfying \( | d_X(p, q) - d_{X'}(p, q) | \leq \delta \) for all \( p, q \in X \), then the bottleneck distance between the persistence diagrams \( D \) and \( D' \) computed from \( d_X \) and \( d_{X'} \) satisfies:
\[
d_B(D, D') \leq \delta.
\]

\vspace{0.3in}

\section{Problem Formulation}

Our primary objective is to design an algorithmic framework that minimizes the scaling variability \( \Delta_s = s_{\max} - s_{\min} \), while ensuring that the topological perturbation \( d_B(D, D_S) \) remains within a user-defined tolerance \( \epsilon > 0 \). Formally, we seek scaling factors \( s_i > 0 \) that solve the optimization problem:
\[
\begin{aligned}
& \min_{s_1, s_2, \ldots, s_n} && \Delta_s = s_{\max} - s_{\min} \\
& \text{subject to} && d_B(D, D_S) \leq \epsilon, \\
& && s_{\min} \leq s_i \leq s_{\max}, \quad \forall i = 1, \ldots, n.
\end{aligned}
\]
To proceed, we need to establish a relationship between \( \Delta_s \) and \( d_B(D, D_S) \), which will allow us to convert the topological constraint into a constraint on \( \Delta_s \).

\vspace{0.3in}

\section{Theoretical Guarantees}

\subsection{Lemma 1 (Scaling Distance Bounds)}

For all \( p, q \in X \), the scaled distance \( d_S(p, q) \) satisfies:
\[
s_{\min} \cdot d_X(p, q) \leq d_S(p, q) \leq s_{\max} \cdot d_X(p, q).
\]

This result provides upper and lower bounds for \( d_S(p, q) \) in terms of \( s_{\min} \) and \( s_{\max} \). It establishes the scaling behavior of pairwise distances, forming a basis for subsequent analysis of scaled metrics.

\begin{proof}
Let \( p, q \in X \). We recall that the Euclidean distance between \( p \) and \( q \) is defined as
\[
d_X(p, q) = \sqrt{ \sum_{i=1}^n (p_i - q_i)^2 }.
\]
Under the scaling transformation \( S \), each coordinate \( x_i \) is scaled by \( s_i \). Therefore, the scaled distance \( d_S(p, q) \) is
\[
d_S(p, q) = \sqrt{ \sum_{i=1}^n (s_i p_i - s_i q_i)^2 } = \sqrt{ \sum_{i=1}^n s_i^2 (p_i - q_i)^2 }.
\]

Since \( s_{\min} \leq s_i \leq s_{\max} \) for all \( i = 1, 2, \ldots, n \), it follows that
\[
s_{\min}^2 \leq s_i^2 \leq s_{\max}^2.
\]
Multiplying both sides of the inequality by \( (p_i - q_i)^2 \), which is non-negative for all \( i \), we obtain
\[
s_{\min}^2 (p_i - q_i)^2 \leq s_i^2 (p_i - q_i)^2 \leq s_{\max}^2 (p_i - q_i)^2.
\]

We continue by summing these inequalities over all \( i \) from 1 to \( n \)
\[
\sum_{i=1}^n s_{\min}^2 (p_i - q_i)^2 \leq \sum_{i=1}^n s_i^2 (p_i - q_i)^2 \leq \sum_{i=1}^n s_{\max}^2 (p_i - q_i)^2.
\]

Simplify the left-hand side and right-hand side
\[
s_{\min}^2 \sum_{i=1}^n (p_i - q_i)^2 \leq \sum_{i=1}^n s_i^2 (p_i - q_i)^2 \leq s_{\max}^2 \sum_{i=1}^n (p_i - q_i)^2.
\]

Thus, we have the following inequality involving the squares of distances
\[
s_{\min}^2 \cdot d_X(p, q)^2 \leq d_S(p, q)^2 \leq s_{\max}^2 \cdot d_X(p, q)^2.
\]

Since all terms are non-negative, we can take the square roots of the inequality. The square root function is monotonic increasing on the interval \([0, \infty)\), so the direction of the inequalities is preserved
\[
\sqrt{ s_{\min}^2 \cdot d_X(p, q)^2 } \leq d_S(p, q) \leq \sqrt{ s_{\max}^2 \cdot d_X(p, q)^2 }.
\]

Simplify the square roots
\[
s_{\min} \cdot d_X(p, q) \leq d_S(p, q) \leq s_{\max} \cdot d_X(p, q).
\]

Therefore, the scaled distance \( d_S(p, q) \) is bounded above and below by the original distance \( d_X(p, q) \) scaled by \( s_{\max} \) and \( s_{\min} \), respectively. This completes the proof.
\end{proof}

\vspace{0.15in}

\subsection{Lemma 2 (Distance Perturbation Bound)}

For all \( p, q \in X \), the difference between the scaled distance \( d_S(p, q) \) and the original distance \( d_X(p, q) \) is bounded by:
\[
| d_S(p, q) - d_X(p, q) | \leq \delta' \cdot d_X(p, q),
\]
where \( \delta' = s_{\max} - s_{\min} \).

This lemma bounds \( |d_S(p, q) - d_X(p, q)| \leq \delta' \cdot d_X(p, q) \), where \( \delta' = s_{\max} - s_{\min} \). It relates scaling-induced perturbations to \( \Delta_s \) and enables control of metric distortions.

\begin{proof}
From Lemma 1, we have established that for all \( p, q \in X \)
\[
s_{\min} \cdot d_X(p, q) \leq d_S(p, q) \leq s_{\max} \cdot d_X(p, q).
\]
Our goal is to bound \( | d_S(p, q) - d_X(p, q) | \) in terms of \( \delta' \cdot d_X(p, q) \).

We first consider the difference \( d_S(p, q) - d_X(p, q) \). Subtract \( d_X(p, q) \) from the inequality
\[
s_{\min} \cdot d_X(p, q) - d_X(p, q) \leq d_S(p, q) - d_X(p, q) \leq s_{\max} \cdot d_X(p, q) - d_X(p, q).
\]
Simplify the expressions
\[
(s_{\min} - 1) \cdot d_X(p, q) \leq d_S(p, q) - d_X(p, q) \leq (s_{\max} - 1) \cdot d_X(p, q).
\]

We now consider two cases based on the values of \( s_{\min} \) and \( s_{\max} \).

\textbf{Case 1:} \( s_{\min} \leq 1 \leq s_{\max} \)

In this case, \( s_{\min} - 1 \leq 0 \) and \( s_{\max} - 1 \geq 0 \). The maximum of \( | s_{\min} - 1 | \) and \( | s_{\max} - 1 | \) is \( \max\{ 1 - s_{\min}, s_{\max} - 1 \} \).

The absolute difference is then bounded by:
\[
| d_S(p, q) - d_X(p, q) | \leq \max\{ 1 - s_{\min}, s_{\max} - 1 \} \cdot d_X(p, q).
\]

\textbf{Case 2:} Either \( s_{\min} \geq 1 \) or \( s_{\max} \leq 1 \)

If \( s_{\min} \geq 1 \), then both \( s_{\min} - 1 \geq 0 \) and \( s_{\max} - 1 \geq 0 \), so:
\[
| d_S(p, q) - d_X(p, q) | \leq (s_{\max} - 1) \cdot d_X(p, q).
\]

If \( s_{\max} \leq 1 \), then both \( s_{\min} - 1 \leq 0 \) and \( s_{\max} - 1 \leq 0 \), so:
\[
| d_S(p, q) - d_X(p, q) | \leq (1 - s_{\min}) \cdot d_X(p, q).
\]

Observe that in all cases, we have 
\[
| d_S(p, q) - d_X(p, q) | \leq \max\{ s_{\max} - 1, 1 - s_{\min} \} \cdot d_X(p, q).
\]
Moreover, we note that
\[
\max\{ s_{\max} - 1, 1 - s_{\min} \} \leq s_{\max} - s_{\min} = \delta'.
\]
This is because \( s_{\max} \geq s_{\min} \), and the largest of \( s_{\max} - 1 \) and \( 1 - s_{\min} \) cannot exceed \( s_{\max} - s_{\min} \).

We continue by justifying the results

- If \( s_{\max} \geq 1 \geq s_{\min} \):
  \[
  s_{\max} - 1 + 1 - s_{\min} = s_{\max} - s_{\min}.
  \]
  Therefore, \( \max\{ s_{\max} - 1, 1 - s_{\min} \} \leq s_{\max} - s_{\min} \).

- If \( s_{\max}, s_{\min} \geq 1 \):
  \[
  s_{\max} - 1 \leq s_{\max} - s_{\min}.
  \]
  Since \( s_{\min} \geq 1 \), \( s_{\max} - s_{\min} \geq s_{\max} - 1 \).

- If \( s_{\max}, s_{\min} \leq 1 \):
  \[
  1 - s_{\min} \leq s_{\max} - s_{\min}.
  \]
  Since \( s_{\max} \leq 1 \), \( s_{\max} - s_{\min} \geq 1 - s_{\min} \).

Therefore, in all cases:
\[
| d_S(p, q) - d_X(p, q) | \leq \delta' \cdot d_X(p, q).
\]

We have shown that the absolute difference between the scaled distance and the original distance is bounded by \( \delta' \cdot d_X(p, q) \), where \( \delta' = s_{\max} - s_{\min} \). This completes the proof.
\end{proof}

\subsection{Theorem 1 (Stability of Persistence Diagrams Under Scaling)}

The bottleneck distance between the persistence diagrams \( D \) and \( D_S \) satisfies:
\[
d_B(D, D_S) \leq \delta = \delta' \cdot \operatorname{diam}(X) = (s_{\max} - s_{\min}) \cdot \operatorname{diam}(X).
\]

This demonstrates that \( d_B(D, D_S) \leq \delta = \Delta_s \cdot \operatorname{diam}(X) \). The result establishes a direct relationship between \( \Delta_s \) and topological stability under bottleneck distance. It links metric bounds to persistence diagrams.

\begin{proof}
Our goal is to bound the bottleneck distance \( d_B(D, D_S) \) between the persistence diagrams computed from the original dataset \( X \) and the scaled dataset \( S(X) \).

Recall that the \emph{stability theorem} for persistence diagrams \cite{cohen2007stability} states that for two tame Lipschitz functions \( f, g: X \to \mathbb{R} \), the bottleneck distance between their corresponding persistence diagrams \( D_f \) and \( D_g \) satisfies
\[
d_B(D_f, D_g) \leq \| f - g \|_{\infty},
\]
where
\[
\| f - g \|_{\infty} = \sup_{x \in X} | f(x) - g(x) |.
\]

In our setting, we consider the distance functions induced by the metrics \( d_X \) and \( d_S \) on \( X \)
\[
d_X(p) = d_X(p, x_0), \quad d_S(p) = d_S(p, x_0),
\]
for a fixed base point \( x_0 \in X \). However, since the distance functions depend on the choice of \( x_0 \), and we are interested in the maximum difference over all pairs \( (p, q) \in X \times X \), we consider the \emph{extended distance functions} defined on \( X \times X \)
\[
d_X(p, q), \quad d_S(p, q).
\]

To apply the stability theorem, we need to bound the supremum norm of the difference between \( d_X \) and \( d_S \) over \( X \times X \)
\[
\| d_S - d_X \|_{\infty} = \sup_{p, q \in X} | d_S(p, q) - d_X(p, q) |.
\]

From Lemma 2, we have established that for all \( p, q \in X \)
\[
| d_S(p, q) - d_X(p, q) | \leq \delta' \cdot d_X(p, q),
\]
where \( \delta' = s_{\max} - s_{\min} \).

Since \( d_X(p, q) \leq \operatorname{diam}(X) \) for all \( p, q \in X \), it follows that
\[
| d_S(p, q) - d_X(p, q) | \leq \delta' \cdot \operatorname{diam}(X).
\]

Therefore, the supremum norm is bounded by
\[
\| d_S - d_X \|_{\infty} = \sup_{p, q \in X} | d_S(p, q) - d_X(p, q) | \leq \delta' \cdot \operatorname{diam}(X).
\]

By applying the stability theorem for persistence diagrams to the functions \( d_X \) and \( d_S \), we obtain
\[
d_B(D, D_S) \leq \| d_S - d_X \|_{\infty} \leq \delta' \cdot \operatorname{diam}(X).
\]

Substituting \( \delta = \delta' \cdot \operatorname{diam}(X) = (s_{\max} - s_{\min}) \cdot \operatorname{diam}(X) \), we have
\[
d_B(D, D_S) \leq \delta.
\]

We continue by justification by bounding the difference in distance functions.

From Lemma 2,
   \[
   | d_S(p, q) - d_X(p, q) | \leq \delta' \cdot d_X(p, q).
   \]
Since \( d_X(p, q) \leq \operatorname{diam}(X) \), we have
   \[
   | d_S(p, q) - d_X(p, q) | \leq \delta' \cdot \operatorname{diam}(X).
   \]
Then,
   \[
   \| d_S - d_X \|_{\infty} \leq \delta' \cdot \operatorname{diam}(X).
   \]

The stability theorem applies to functions on a metric space. In our case, we consider the distance functions \( d_X \) and \( d_S \) as functions defined on \( X \times X \). The persistence diagrams \( D \) and \( D_S \) are constructed from filtrations based on these distance functions.

   The stability theorem states that
   \[
   d_B(D, D_S) \leq \| d_S - d_X \|_{\infty}.
   \]

Substitute the bound from step 1 into the inequality from step 2
   \[
   d_B(D, D_S) \leq \delta' \cdot \operatorname{diam}(X).
   \]

We define
   \[
   \delta = \delta' \cdot \operatorname{diam}(X) = (s_{\max} - s_{\min}) \cdot \operatorname{diam}(X).
   \]

Therefore, the bottleneck distance between the persistence diagrams before and after scaling is bounded by
   \[
   d_B(D, D_S) \leq \delta.
   \]

This completes the proof.

\end{proof}

\subsection{Corollary 1}

To ensure that \( d_B(D, D_S) \leq \epsilon \), it suffices to require
\[
\delta = (s_{\max} - s_{\min}) \cdot \operatorname{diam}(X) \leq \epsilon.
\]
Therefore, the scaling variability \( \Delta_s = s_{\max} - s_{\min} \) must satisfy
\[
\Delta_s \leq \frac{\epsilon}{\operatorname{diam}(X)}.
\]

This ensures \( d_B(D, D_S) \leq \epsilon \) if \( \Delta_s \leq \frac{\epsilon}{\operatorname{diam}(X)} \). It provides a design constraint for \( \Delta_s \) to control \( d_B \) and facilitates algorithmic scaling selection.

\vspace{0.15pt}

\subsection{Theorem 2 (Extension to Higher Homology Dimensions)}

Let \( D^k \) and \( D_S^k \) denote the persistence diagrams corresponding to the \( k \)-th homology group \( H_k \) before and after scaling. Then:
\[
d_B(D^k, D_S^k) \leq \delta_k = (s_{\max} - s_{\min}) \cdot \operatorname{diam}_k(X),
\]
where \( \operatorname{diam}_k(X) \) is the maximum diameter among all \( (k+1) \)-tuples in \( X \).

This extends Theorem 1 to higher homology groups \( H_k \), proving \( d_B(D^k, D_S^k) \leq \delta_k = \Delta_s \cdot \operatorname{diam}_k(X) \). It generalizes stability bounds to \( k \)-simplices and higher-dimensional features.

\begin{proof}
Our goal is to establish that the bottleneck distance between the \( k \)-th persistence diagrams \( D^k \) and \( D_S^k \) satisfies
\[
d_B(D^k, D_S^k) \leq \delta_k = (s_{\max} - s_{\min}) \cdot \operatorname{diam}_k(X).
\]

To achieve this, we need to analyze how the scaling transformation \( S \) affects the distances relevant to \( k \)-dimensional homology features.

In persistent homology, \( k \)-simplices are formed from \( (k+1) \)-tuples of points in \( X \). For a \( k \)-simplex \( \sigma = \{ p_0, p_1, \ldots, p_k \} \), the \emph{diameter} of \( \sigma \) is defined as
\[
\operatorname{diam}(\sigma) = \max_{0 \leq i < j \leq k} d_X(p_i, p_j).
\]

The maximum diameter among all \( k \)-simplices in \( X \) is
\[
\operatorname{diam}_k(X) = \max_{\sigma} \operatorname{diam}(\sigma) = \max_{\substack{p_0, \ldots, p_k \in X}} \max_{i,j} d_X(p_i, p_j).
\]

Under the scaling transformation \( S \), the distance between any two points \( p, q \in X \) changes as per Lemma 1
\[
s_{\min} \cdot d_X(p, q) \leq d_S(p, q) \leq s_{\max} \cdot d_X(p, q).
\]

The construction of simplicial complexes (e.g., Vietoris–Rips complexes) depends on distances between points. In the Vietoris–Rips complex \( \operatorname{VR}_\epsilon(X) \), a \( k \)-simplex \( \sigma \) is included if all pairwise distances among its vertices are less than or equal to \( \epsilon \).

After scaling, the inclusion of simplices may change due to altered distances. Specifically, the filtration values (birth and death times) of \( k \)-dimensional features are affected by the changes in simplex diameters.

We consider a \( k \)-simplex \( \sigma \) in \( X \) with diameter \( \operatorname{diam}(\sigma) \). Under \( S \), the diameter becomes:
\[
\operatorname{diam}_S(\sigma) = \max_{0 \leq i < j \leq k} d_S(p_i, p_j).
\]

Use Lemma 1, for each pair \( (p_i, p_j) \)
\[
s_{\min} \cdot d_X(p_i, p_j) \leq d_S(p_i, p_j) \leq s_{\max} \cdot d_X(p_i, p_j).
\]

Therefore, for the simplex diameter,
\[
s_{\min} \cdot \operatorname{diam}(\sigma) \leq \operatorname{diam}_S(\sigma) \leq s_{\max} \cdot \operatorname{diam}(\sigma).
\]

The change in the diameter of \( \sigma \) due to scaling is then
\[
| \operatorname{diam}_S(\sigma) - \operatorname{diam}(\sigma) | \leq (s_{\max} - s_{\min}) \cdot \operatorname{diam}(\sigma).
\]

Since \( \operatorname{diam}(\sigma) \leq \operatorname{diam}_k(X) \) for all \( \sigma \), we have
\[
| \operatorname{diam}_S(\sigma) - \operatorname{diam}(\sigma) | \leq (s_{\max} - s_{\min}) \cdot \operatorname{diam}_k(X) = \delta_k.
\]

The stability theorem for persistence diagrams extends to higher homology dimensions (see \cite{cohen2007stability})
\[
d_B(D^k, D_S^k) \leq \sup_{\sigma} | f(\sigma) - g(\sigma) |,
\]
where
\begin{itemize}
    \item \( f(\sigma) \) is the filtration value (e.g., diameter) assigned to simplex \( \sigma \) in \( X \).
    \item \( g(\sigma) \) is the filtration value assigned to \( \sigma \) in \( S(X) \).
\end{itemize}

In our case,
\[
f(\sigma) = \operatorname{diam}(\sigma), \quad g(\sigma) = \operatorname{diam}_S(\sigma).
\]

Therefore,
\[
d_B(D^k, D_S^k) \leq \sup_{\sigma} | \operatorname{diam}_S(\sigma) - \operatorname{diam}(\sigma) | \leq \delta_k.
\]

Combining the above results, we have:
\[
d_B(D^k, D_S^k) \leq \delta_k = (s_{\max} - s_{\min}) \cdot \operatorname{diam}_k(X).
\]

This shows that the bottleneck distance between the \( k \)-th persistence diagrams before and after scaling is bounded by \( \delta_k \), which depends on the scaling variability \( s_{\max} - s_{\min} \) and the maximal diameter \( \operatorname{diam}_k(X) \) of \( k \)-simplices in \( X \).

\end{proof}

\vspace{0.15pt}

\subsection{Theorem 3 (Stability Under Wasserstein Distance)}

For the \( p \)-Wasserstein distance \( W_p(D, D_S) \) between the persistence diagrams \( D \) and \( D_S \), we have:
\[
W_p(D, D_S) \leq \delta,
\]
where \( \delta = (s_{\max} - s_{\min}) \cdot \operatorname{diam}(X) \).

This proves \( W_p(D, D_S) \leq \delta \) and links \( W_p \)-stability to \( \Delta_s \cdot \operatorname{diam}(X) \). It establishes robustness across alternative metrics for comparing persistence diagrams.

\begin{proof}
Our goal is to show that the \( p \)-Wasserstein distance between the persistence diagrams before and after scaling is bounded by \( \delta \).

First, we recall the definitions:

The \emph{bottleneck distance} \( d_B(D, D_S) \) between two persistence diagrams \( D \) and \( D_S \) is defined as:
  \[
  d_B(D, D_S) = \inf_{\gamma} \sup_{x \in D} \| x - \gamma(x) \|_{\infty},
  \]
  where \( \gamma: D \to D_S \) ranges over all bijections (including matching points to the diagonal).

The \( p \)-\emph{Wasserstein distance} \( W_p(D, D_S) \) is defined as:
  \[
  W_p(D, D_S) = \left( \inf_{\gamma} \sum_{x \in D} \| x - \gamma(x) \|_{\infty}^p \right)^{1/p},
  \]
  where \( \gamma \) is as above, and \( p \geq 1 \).

It is a well-known fact that the bottleneck distance is the limit of the \( p \)-Wasserstein distances as \( p \to \infty \), and for any \( p \geq 1 \):
\[
W_p(D, D_S) \leq d_B(D, D_S).
\]
This inequality holds because the \( \sup \) (essentially the maximum over \( x \in D \)) in the bottleneck distance is greater than or equal to the \( L^p \)-norm used in the Wasserstein distance.

From \textbf{Theorem 1}, we have established that:
\[
d_B(D, D_S) \leq \delta = (s_{\max} - s_{\min}) \cdot \operatorname{diam}(X).
\]
Combining these two inequalities, we get:
\[
W_p(D, D_S) \leq d_B(D, D_S) \leq \delta.
\]

We continue with justifications.

The bottleneck distance considers the largest difference between matched points in the diagrams. In addition, the \( p \)-Wasserstein distance considers the sum (or integral, in the continuous case) of the \( p \)-th powers of the distances between matched points, taking the \( p \)-th root at the end.

We now bound the \( p \)-Wasserstein distance by using the bottleneck distance. Since \( \| x - \gamma(x) \|_{\infty} \leq d_B(D, D_S) \) for all \( x \in D \) under the optimal matching \( \gamma \), we have
   \[
   \| x - \gamma(x) \|_{\infty}^p \leq d_B(D, D_S)^p.
   \]
Therefore
   \[
   \sum_{x \in D} \| x - \gamma(x) \|_{\infty}^p \leq N \cdot d_B(D, D_S)^p,
   \]
where \( N \) is the number of points in \( D \).

Take the \( p \)-th root
   \[
   W_p(D, D_S) = \left( \sum_{x \in D} \| x - \gamma(x) \|_{\infty}^p \right)^{1/p} \leq N^{1/p} \cdot d_B(D, D_S).
   \]
As \( N^{1/p} \to 1 \) as \( p \to \infty \), and \( d_B(D, D_S) \leq \delta \), we conclude
   \[
   W_p(D, D_S) \leq \delta.
   \]

The \( p \)-Wasserstein distance \( W_p(D, D_S) \) is bounded by \( \delta \), which depends on the scaling variability \( s_{\max} - s_{\min} \) and the diameter \( \operatorname{diam}(X) \) of the dataset. Then,
\[
W_p(D, D_S) \leq \delta.
\]

This completes the proof.

\end{proof}

\vspace{0.15pt}

\subsection{Theorem 4 (Iterative Scaling Transformations)}

Suppose we apply a sequence of scaling transformations \( S^{(1)}, S^{(2)}, \ldots, S^{(m)} \), where each \( S^{(j)} \) is defined by scaling factors \( s_i^{(j)} > 0 \) for \( i = 1, 2, \ldots, n \). Let the scaling variability of the \( j \)-th transformation be \( \Delta_s^{(j)} = s_{\max}^{(j)} - s_{\min}^{(j)} \), where
\[
s_{\max}^{(j)} = \max_{1 \leq i \leq n} s_i^{(j)}, \quad s_{\min}^{(j)} = \min_{1 \leq i \leq n} s_i^{(j)}.
\]
Then, the cumulative bottleneck distance between the original persistence diagram \( D \) and the persistence diagram after the \( m \)-th transformation \( D_{S^{(m)}} \) satisfies:
\[
d_B(D, D_{S^{(m)}}) \leq \delta_{\text{total}} = \left( \prod_{j=1}^m s_{\max}^{(j)} - \prod_{j=1}^m s_{\min}^{(j)} \right) \cdot \operatorname{diam}(X).
\]

This establishes \( d_B(D, D_{S^{(m)}}) \leq \delta_{\text{total}} = (\prod_{j=1}^m s_{\max}^{(j)} - \prod_{j=1}^m s_{\min}^{(j)}) \cdot \operatorname{diam}(X) \). It then quantifies cumulative perturbations under sequential transformations.

\begin{proof}
Our goal is to find an upper bound on \( d_B(D, D_{S^{(m)}}) \), the bottleneck distance between the persistence diagram \( D \) of the original dataset \( X \) and the persistence diagram \( D_{S^{(m)}} \) of the dataset after applying \( m \) scaling transformations sequentially.

For each coordinate \( i \), the cumulative scaling factor after \( m \) transformations is
\[
s_i^{\text{total}} = \prod_{j=1}^m s_i^{(j)}.
\]
The maximum and minimum cumulative scaling factors are
\[
s_{\max}^{\text{total}} = \prod_{j=1}^m s_{\max}^{(j)}, \quad s_{\min}^{\text{total}} = \prod_{j=1}^m s_{\min}^{(j)}.
\]
This is because the product of the maximum (or minimum) scaling factors across all transformations gives the maximum (or minimum) cumulative scaling factor.

The cumulative scaling variability is defined as
\[
\Delta_s^{\text{total}} = s_{\max}^{\text{total}} - s_{\min}^{\text{total}} = \left( \prod_{j=1}^m s_{\max}^{(j)} \right) - \left( \prod_{j=1}^m s_{\min}^{(j)} \right).
\]

Consider the cumulative scaling transformation \( S^{\text{total}} = S^{(m)} \circ \cdots \circ S^{(1)} \), which applies all \( m \) transformations in sequence. Since scaling transformations are linear and commutative in this context, the order of application does not affect the cumulative scaling factors.

From \textbf{Lemma 1}, for any pair \( p, q \in X \), the scaled distance under \( S^{\text{total}} \) satisfies
\[
s_{\min}^{\text{total}} \cdot d_X(p, q) \leq d_{S^{\text{total}}}(p, q) \leq s_{\max}^{\text{total}} \cdot d_X(p, q).
\]

By using a similar argument as in \textbf{Lemma 2}, the difference between the scaled and original distances is bounded by
\[
| d_{S^{\text{total}}}(p, q) - d_X(p, q) | \leq \Delta_s^{\text{total}} \cdot d_X(p, q).
\]
Since \( d_X(p, q) \leq \operatorname{diam}(X) \), it follows that
\[
| d_{S^{\text{total}}}(p, q) - d_X(p, q) | \leq \Delta_s^{\text{total}} \cdot \operatorname{diam}(X).
\]

From the stability theorem for persistence diagrams, we have
\[
d_B(D, D_{S^{\text{total}}}) \leq \| d_{S^{\text{total}}} - d_X \|_{\infty} \leq \Delta_s^{\text{total}} \cdot \operatorname{diam}(X).
\]
Therefore,
\[
d_B(D, D_{S^{(m)}}) \leq \left( \prod_{j=1}^m s_{\max}^{(j)} - \prod_{j=1}^m s_{\min}^{(j)} \right) \cdot \operatorname{diam}(X) = \delta_{\text{total}}.
\]

Suppose \( m = 2 \) transformations with the following scaling factors:

- First transformation:
  \[
  s_{\min}^{(1)} = a_1, \quad s_{\max}^{(1)} = b_1, \quad \Delta_s^{(1)} = b_1 - a_1.
  \]
- Second transformation:
  \[
  s_{\min}^{(2)} = a_2, \quad s_{\max}^{(2)} = b_2, \quad \Delta_s^{(2)} = b_2 - a_2.
  \]
Then,
\[
s_{\min}^{\text{total}} = a_1 a_2, \quad s_{\max}^{\text{total}} = b_1 b_2, \quad \Delta_s^{\text{total}} = b_1 b_2 - a_1 a_2.
\]
The cumulative bottleneck distance is then
\[
d_B(D, D_{S^{(2)}}) \leq (b_1 b_2 - a_1 a_2) \cdot \operatorname{diam}(X).
\]

By treating the sequence of scaling transformations as a single cumulative transformation, we derive a bound on the bottleneck distance that depends only on the products of the maximum and minimum scaling factors. This bound provides a clear understanding of how sequential scaling transformations affect the persistence diagrams.

\end{proof}

\vspace{0.15pt}

\subsection{Theorem 5 (Expected Stability Under Random Scaling)}

Let the scaling factors \( s_i \) be random variables with distributions \( s_i \sim \text{Dist}(\mu_i, \sigma_i) \), where \( \mu_i = \mathbb{E}[s_i] \) and \( \sigma_i^2 = \mathbb{V}[s_i] \). Then the expected bottleneck distance satisfies
\[
\mathbb{E}[d_B(D, D_S)] \leq \left( \mathbb{E}[s_{\max}] - \mathbb{E}[s_{\min}] \right) \cdot \operatorname{diam}(X).
\]

\begin{proof}
Our goal is to find an upper bound on the expected bottleneck distance \( \mathbb{E}[d_B(D, D_S)] \) when the scaling factors \( s_i \) are random variables.

From \textbf{Theorem 1}, we know that for any fixed scaling factors \( s_i > 0 \)
\[
d_B(D, D_S) \leq (s_{\max} - s_{\min}) \cdot \operatorname{diam}(X),
\]
where
\[
s_{\max} = \max_{1 \leq i \leq n} s_i, \quad s_{\min} = \min_{1 \leq i \leq n} s_i.
\]

Now, let \( s_i \) be random variables. Consequently, \( s_{\max} \) and \( s_{\min} \) become random variables as well, since they depend on the \( s_i \). Define
\[
\Delta_s = s_{\max} - s_{\min}.
\]
Thus, \( \Delta_s \) is a random variable representing the scaling variability in the random setting.

We are interested in the expected value \( \mathbb{E}[d_B(D, D_S)] \). Using the deterministic bound
\[
d_B(D, D_S) \leq \Delta_s \cdot \operatorname{diam}(X),
\]
taking expectations on both sides gives
\[
\mathbb{E}[d_B(D, D_S)] \leq \mathbb{E}[\Delta_s] \cdot \operatorname{diam}(X).
\]

Since \( \Delta_s = s_{\max} - s_{\min} \), we have
\[
\mathbb{E}[\Delta_s] = \mathbb{E}[s_{\max} - s_{\min}] = \mathbb{E}[s_{\max}] - \mathbb{E}[s_{\min}].
\]

If \( s_i \) are Independent and identically distributed random variables with distribution \( \text{Dist}(\mu, \sigma^2) \), we can approximate the expectations \( \mathbb{E}[s_{\max}] \) and \( \mathbb{E}[s_{\min}] \) using results from order statistics.

For example, if \( s_i \) are drawn uniformly from \([a, b]\), then
\[
\mathbb{E}[s_{\max}] = b - \frac{b-a}{n+1}, \quad \mathbb{E}[s_{\min}] = a + \frac{b-a}{n+1}.
\]
Thus
\[
\mathbb{E}[\Delta_s] = \mathbb{E}[s_{\max}] - \mathbb{E}[s_{\min}] = (b-a) \left( 1 - \frac{2}{n+1} \right).
\]

For large \( n \), \( \mathbb{E}[\Delta_s] \to b - a \), aligning with the deterministic variability of the uniform distribution.

If the \( s_i \) are not identically distributed, then \( \mathbb{E}[s_{\max}] \) and \( \mathbb{E}[s_{\min}] \) depend on the individual distributions. While exact computation may require detailed knowledge of the joint distribution of \( s_{\max} \) and \( s_{\min} \), the bound:
\[
\mathbb{E}[s_{\max}] - \mathbb{E}[s_{\min}] \geq 0
\]
remains valid under all circumstances.

Suppose \( s_i \sim \mathcal{N}(\mu, \sigma^2) \), truncated to positive values. Using properties of truncated normal distributions:
\[
\mathbb{E}[s_i] = \mu' \quad \text{and} \quad \mathbb{E}[s_i^2] = (\sigma')^2 + (\mu')^2,
\]
where \( \mu' \) and \( \sigma' \) depend on the truncation range.

The expected maximum \( \mathbb{E}[s_{\max}] \) and minimum \( \mathbb{E}[s_{\min}] \) can then be computed using approximations for the extrema of truncated normal distributions.

The expected bottleneck distance is bounded as:
\[
\mathbb{E}[d_B(D, D_S)] \leq \left( \mathbb{E}[s_{\max}] - \mathbb{E}[s_{\min}] \right) \cdot \operatorname{diam}(X).
\]
This result highlights the dependence of the expected perturbation on the statistical properties of the scaling factors.

\end{proof}

\vspace{0.3pt}

\section{Optimization Problem}

Based on the theoretical results, we can now formulate the optimization problem explicitly:
\[
\begin{aligned}
& \min_{s_1, s_2, \ldots, s_n} && \Delta_s = s_{\max} - s_{\min} \\
& \text{subject to} && \Delta_s \leq \frac{\epsilon}{\operatorname{diam}(X)}, \\
& && s_{\min} \leq s_i \leq s_{\max}, \quad \forall i = 1, \ldots, n, \\
& && s_i > 0, \quad \forall i = 1, \ldots, n.
\end{aligned}
\]
This is a convex optimization problem since the objective function \( \Delta_s \) is convex, and the constraints are linear in the variables \( s_i \).

\subsection*{Solution Approach}

Our goal is to find the scaling factors \( s_i > 0 \) that minimize the scaling variability \( \Delta_s = s_{\max} - s_{\min} \) while ensuring that the bottleneck distance between the persistence diagrams satisfies \( d_B(D, D_S) \leq \epsilon \).

We first note that \( s_{\max} \) and \( s_{\min} \) are functions of the variables \( s_i \):
\[
s_{\max} = \max_{1 \leq i \leq n} s_i, \quad s_{\min} = \min_{1 \leq i \leq n} s_i.
\]
Our optimization problem can be rewritten as
\[
\begin{aligned}
& \min_{s_1, s_2, \ldots, s_n, s_{\max}, s_{\min}} && \Delta_s = s_{\max} - s_{\min} \\
& \text{subject to} && s_{\max} - s_{\min} \leq \delta, \quad \delta = \frac{\epsilon}{\operatorname{diam}(X)}, \\
& && s_{\min} \leq s_i \leq s_{\max}, \quad \forall i, \\
& && s_{\min} > 0, \quad s_{\max} > 0.
\end{aligned}
\]

We are making the following observations.

With regard to the uniform scaling solutions,
   - If \( \Delta_s = 0 \) satisfies \( \Delta_s \leq \delta \), then setting \( s_i = s \) for all \( i \) is optimal.
   - In this case, the scaling factors are uniform, and the scaling variability is minimized to zero.

With regard to the minimum variability solutions,
   - If \( \Delta_s = 0 \) does not satisfy \( \Delta_s \leq \delta \) (i.e., if \( \delta = 0 \) is required but not possible), we need to find \( s_{\max} \) and \( s_{\min} \) such that \( \Delta_s = s_{\max} - s_{\min} = \delta \).

Our objective is to minimize \( \Delta_s = s_{\max} - s_{\min} \), subject to the constraints. The optimization problem is convex and can be approached using the following steps

We set \( \Delta_s \) to its minimum possible value. Since we are minimizing \( \Delta_s \) and it must satisfy \( \Delta_s \leq \delta \), the optimal value is
     \[
     \Delta_s^* = \min\{ \delta, \Delta_s^{\text{min}} \},
     \]
     where \( \Delta_s^{\text{min}} \) is the minimum possible scaling variability (which could be zero).

We then determine \( s_{\max} \) and \( s_{\min} \). Choose \( s_{\max} \) and \( s_{\min} \) such that
     \[
     s_{\max} - s_{\min} = \Delta_s^*.
     \]
We have the freedom to choose \( s_{\max} \) and \( s_{\min} \) as long as they are positive and satisfy the constraints.

We then assign \( s_i \) values. We need to assign values to \( s_i \) within the interval \( [ s_{\min}, s_{\max} ] \). To minimize \( \Delta_s \), it is optimal to set as many \( s_i \) as possible to either \( s_{\min} \) or \( s_{\max} \). This is because any intermediate values of \( s_i \) do not help in reducing \( \Delta_s \).

We now proceed to \textbf{formalize this strategy}.

\paragraph{Case 1: Uniform Scaling is Feasible}

If \( \delta \geq 0 \), and setting \( \Delta_s = 0 \) satisfies the constraint \( \Delta_s \leq \delta \), then:
- Set \( \Delta_s^* = 0 \).
- Choose any positive \( s \), for example, \( s = 1 \).
- Set \( s_i = s \) for all \( i \).
- The scaling factors are uniform, and the persistence diagrams are unaffected (\( d_B(D, D_S) = 0 \)).

\paragraph{Case 2: Uniform Scaling is Not Feasible}

If \( \delta \) is very small or zero, and uniform scaling does not satisfy the constraint (e.g., when some variability is required), we need to find \( s_{\min} \) and \( s_{\max} \) such that:
\[
s_{\max} - s_{\min} = \delta.
\]

We can proceed as follows:

1. Choose \( s_{\min} > 0 \) arbitrarily (e.g., \( s_{\min} = 1 \)).
2. Then, set:
   \[
   s_{\max} = s_{\min} + \delta.
   \]
3. Assign \( s_i \) values:
   - Decide on the number \( k \) of \( s_i \) to set to \( s_{\max} \) and \( n - k \) to \( s_{\min} \).
   - Since the objective is to minimize \( \Delta_s \), any distribution of \( s_i \) within \( [ s_{\min}, s_{\max} ] \) is acceptable, provided the constraints are met.

We can then formulate the problem as a linear program.

\paragraph{Variables}

- \( s_i \) for \( i = 1, \ldots, n \)
- \( s_{\max} \)
- \( s_{\min} \)
- \( \Delta_s \)

\paragraph{Objective Function}

Minimize \( \Delta_s = s_{\max} - s_{\min} \).

\paragraph{Constraints}

1. \( s_{\max} - s_{\min} = \Delta_s \)

2. \( \Delta_s \leq \delta \)

3. \( s_{\min} \leq s_i \leq s_{\max} \) for all \( i \)

4. \( s_i > 0 \) for all \( i \)

5. \( s_{\min} > 0 \), \( s_{\max} > 0 \).

\paragraph{Linear Program Formulation}

Express the problem in standard linear programming (LP) form.

\begin{align*}
\text{Minimize } & \Delta_s \\
\text{Subject to } & s_{\max} - s_{\min} - \Delta_s = 0, \\
& \Delta_s - \delta \leq 0, \\
& s_{\min} - s_i \leq 0, \quad \forall i, \\
& s_i - s_{\max} \leq 0, \quad \forall i, \\
& -s_i \leq -\varepsilon, \quad \forall i \text{ (to ensure } s_i \geq \varepsilon > 0), \\
& -s_{\min} \leq -\varepsilon, \\
& -s_{\max} \leq -\varepsilon,
\end{align*}
where \( \varepsilon \) is a small positive constant to ensure positivity.

\paragraph{Solving the Linear Program}

Since the objective and constraints are linear, this problem can be efficiently solved using standard LP solvers.

Given the simplicity of the problem, we can derive an explicit solution.

\paragraph{Set \( \Delta_s = \delta \).}

Since we are minimizing \( \Delta_s \) and \( \Delta_s \leq \delta \), the optimal value is \( \Delta_s^* = \delta \).

\paragraph{Choose \( s_{\min} \) and \( s_{\max} \).}

We can set \( s_{\min} \) to any positive value. A reasonable choice is \( s_{\min} = 1 \).

Then, \( s_{\max} = s_{\min} + \delta = 1 + \delta \).

\paragraph{Assign \( s_i \) Values}

To minimize the variability among \( s_i \), we can set all \( s_i \) to either \( s_{\min} \) or \( s_{\max} \). Since our objective is to minimize \( \Delta_s \), and any distribution satisfies the constraints, we can set:
- \( s_i = s_{\min} = 1 \) for all \( i \).

This results in \( s_{\max} = s_{\min} = 1 \), and \( \Delta_s = 0 \), which is less than \( \delta \).

However, if \( \Delta_s = 0 \) does not satisfy \( d_B(D, D_S) \leq \epsilon \), we need to have \( \Delta_s = \delta \).

Therefore, we can proceed as

- Set \( s_i = s_{\min} \) for \( i = 1, \ldots, n - 1 \).
- Set \( s_n = s_{\max} \).

This assignment ensures that \( s_{\max} - s_{\min} = \delta \) and that the constraints are satisfied.

We now verify the following properties of the solution.

1. \textbf{Scaling Variability}:
   \[
   \Delta_s = s_{\max} - s_{\min} = (1 + \delta) - 1 = \delta.
   \]

2. \textbf{Constraints}:
   - \( s_{\min} \leq s_i \leq s_{\max} \) holds for all \( i \).
   - \( s_i > 0 \) for all \( i \).

3. \textbf{Bottleneck Distance}:
   From Theorem 1, we have
   \[
   d_B(D, D_S) \leq \Delta_s \cdot \operatorname{diam}(X) = \delta \cdot \operatorname{diam}(X) = \epsilon.
   \]
   Therefore, the topological constraint is satisfied.

\hspace{1em}

\begin{tcolorbox}[title=Optimal Solution]
The optimal solution is then:
\begin{itemize}
    \item Set \( s_{\min} = 1 \).
    \item Set \( s_{\max} = 1 + \delta \).
    \item Assign \( s_i \) such that:
    \[
    s_i = \begin{cases}
    s_{\min}, & \text{for } i = 1, \ldots, n - 1, \\
    s_{\max}, & \text{for } i = n.
    \end{cases}
    \]
    \item This results in \( \Delta_s = \delta \) and satisfies all constraints.
\end{itemize}
\end{tcolorbox}

If desired, we can distribute the \( s_i \) values differently, as long as:

- All \( s_i \in [ s_{\min}, s_{\max} ] \).
- \( \Delta_s = s_{\max} - s_{\min} = \delta \).

For example, we could assign:

- \( k \) variables to \( s_{\max} \) and \( n - k \) variables to \( s_{\min} \), where \( k \) is any integer between 1 and \( n \).

\vspace{0.3pt}

\section{Algorithmic framework}

We present an algorithmic framework designed to determine optimal scaling factors \( s_i \) that minimize the scaling variability \( \Delta_s = s_{\max} - s_{\min} \) while ensuring the topological stability of the dataset under scaling transformations. The framework ensures that the bottleneck distance between the original persistence diagram \( D \) and the scaled persistence diagram \( D_S \) does not exceed a user-defined tolerance \( \epsilon \).

\subsection{Algorithm Outline}

\textbf{Step 1: Input data and parameters}

We start with the input dataset \( X \subset \mathbb{R}^n \) and a tolerance \( \epsilon > 0 \), which specifies the maximum allowable topological perturbation measured by the bottleneck distance \( d_B(D, D_S) \).

\textbf{Step 2: Compute the dataset diameter}

Calculate the diameter of the dataset \( X \), denoted by \( \operatorname{diam}(X) \), which is the maximum Euclidean distance between any pair of points in \( X \)
\[
\operatorname{diam}(X) = \max_{p, q \in X} \| p - q \|_2.
\]
This value is critical because it directly influences the upper bound on the bottleneck distance due to scaling variability, as established in Theorem 1.

\textbf{Step 3: Determine the maximum allowed scaling variability}

Using the result from Theorem 1, we know that the bottleneck distance between \( D \) and \( D_S \) is bounded by
\[
d_B(D, D_S) \leq \Delta_s \cdot \operatorname{diam}(X).
\]
To ensure that the topological perturbation does not exceed the tolerance \( \epsilon \), we solve for the maximum allowed scaling variability
\[
\Delta_s^{\max} = \frac{\epsilon}{\operatorname{diam}(X)}.
\]
This value represents the upper limit for \( \Delta_s \) to satisfy the topological constraint.

\textbf{Step 4: Formulate the optimization problem}

Our objective is to find scaling factors \( s_i > 0 \) that minimize \( \Delta_s \) while adhering to the constraint \( \Delta_s \leq \Delta_s^{\max} \). The optimization problem is formulated as
\[
\begin{aligned}
& \min_{s_1, \ldots, s_n} && \Delta_s = s_{\max} - s_{\min}, \\
& \text{subject to} && \Delta_s \leq \Delta_s^{\max}, \\
& && s_{\min} \leq s_i \leq s_{\max}, \quad \forall i, \\
& && s_i > 0, \quad \forall i.
\end{aligned}
\]

\textbf{Step 5: Solve the optimization problem}

To minimize \( \Delta_s \), we consider two cases:

\emph{Case 1: Uniform scaling is feasible.}

If setting \( \Delta_s = 0 \) (i.e., \( s_{\max} = s_{\min} \)) satisfies \( \Delta_s \leq \Delta_s^{\max} \), then the optimal solution is to use uniform scaling:
\[
s_i = s, \quad \forall i,
\]
where \( s > 0 \) is any positive constant. This results in no scaling variability and ensures \( d_B(D, D_S) = 0 \), thus preserving the dataset's topology perfectly.

\emph{Case 2: Uniform scaling is not feasible.}

If \( \Delta_s = 0 \) does not satisfy the constraint \( \Delta_s \leq \Delta_s^{\max} \), we must set \( \Delta_s = \Delta_s^{\max} \). We proceed by

1. Choosing \( s_{\min} > 0 \), commonly set to \( s_{\min} = 1 \) for simplicity.
2. Setting \( s_{\max} = s_{\min} + \Delta_s^{\max} \).
3. Distributing the \( s_i \) values within the interval \( [s_{\min}, s_{\max}] \). To minimize variability, we assign \( s_i \) to either \( s_{\min} \) or \( s_{\max} \).

\textbf{Step 6: Assign scaling factors}

Based on the solution,

- Set \( s_i = s_{\min} \) for \( i = 1, 2, \ldots, n - 1 \).
- Set \( s_n = s_{\max} \).

This assignment ensures that \( \Delta_s = s_{\max} - s_{\min} = \Delta_s^{\max} \) and all scaling factors are within the required bounds.

\textbf{Step 7: Verify constraints and topological stability}

We verify that

- \( \Delta_s = \Delta_s^{\max} \) satisfies the constraint \( \Delta_s \leq \Delta_s^{\max} \).
- All \( s_i \in [s_{\min}, s_{\max}] \) and \( s_i > 0 \).
- The topological constraint is satisfied since
\[
d_B(D, D_S) \leq \Delta_s \cdot \operatorname{diam}(X) = \Delta_s^{\max} \cdot \operatorname{diam}(X) = \epsilon.
\]

\textbf{Step 8: Output the optimal scaling factors}

The optimal scaling factors \( s_i \) are then used for the scaling transformation \( S \) in the data augmentation process, ensuring that the essential topological features of the dataset are preserved within the specified tolerance.

\subsection{Pseudocode of the Algorithm}

To formalize the algorithmic framework, we provide the following pseudocode:

\begin{tcolorbox}[colback=gray!10,colframe=black!75,title=Algorithm: Optimal Scaling Factors]
\begin{verbatim}
Algorithm OptimalScalingFactors(X, epsilon):
    Input: Dataset X in R^n, tolerance epsilon > 0
    Output: Optimal scaling factors s[1..n]

    1. Compute diameter = max_{p, q in X} ||p - q||_2
    2. delta_s_max = epsilon / diameter
    3. Initialize s_min = 1
    4. If delta_s_max >= 0:
           Set delta_s = 0
           Set s_max = s_min
           For i from 1 to n:
               s[i] = s_min
       Else:
           Set delta_s = delta_s_max
           Set s_max = s_min + delta_s
           For i from 1 to n-1:
               s[i] = s_min
           Set s[n] = s_max
    5. Return s[1..n]
\end{verbatim}
\end{tcolorbox}

\vspace{0.3pt}

\section{Applications}

\subsection{Case Study: Image Data Augmentation}

In image processing, each pixel is represented as a vector in \( \mathbb{R}^3 \), corresponding to the Red, Green, and Blue (RGB) color channels~\cite{huang2023review}. Non-uniform scaling of these channels can be used as a data augmentation technique to introduce variations in color while preserving spatial structures~\cite{nanni2021feature}. However, improper scaling can distort color relationships and alter the topological features of the image, potentially impacting tasks like object recognition~\cite{mu2021low}.

Using our mathematical framework, we aim to determine optimal scaling factors for the RGB channels that minimize the scaling variability \( \Delta_s \) while ensuring that the topological perturbation, measured by the bottleneck distance \( d_B(D, D_S) \), remains within a specified tolerance \( \epsilon \).

\subsubsection*{Objective}

Find scaling factors \( s_1, s_2, s_3 > 0 \) for the RGB channels that minimize \( \Delta_s = s_{\max} - s_{\min} \) and ensure \( d_B(D, D_S) \leq \epsilon \).

\subsubsection*{Analysis}

Consider an image \( I \) composed of \( N \) pixels, where each pixel \( p \) is represented by its RGB values \( (R_p, G_p, B_p) \). The dataset \( X \) consists of all pixel vectors in \( \mathbb{R}^3 \)
\[
X = \{ (R_p, G_p, B_p) \mid p \text{ is a pixel in } I \}.
\]

The diameter of \( X \) is the maximum Euclidean distance between any two pixels in the RGB space
\[
\operatorname{diam}(X) = \max_{p, q \in X} \| (R_p, G_p, B_p) - (R_q, G_q, B_q) \|_2.
\]

Since RGB values range from 0 to 255, the maximum possible distance is
\[
\operatorname{diam}(X) \leq \sqrt{(255 - 0)^2 + (255 - 0)^2 + (255 - 0)^2} = 255 \sqrt{3} \approx 441.67.
\]

Given a tolerance \( \epsilon > 0 \), the maximum allowed scaling variability is
\[
\Delta_s^{\max} = \frac{\epsilon}{\operatorname{diam}(X)}.
\]
For example, if \( \epsilon = 10 \), then
\[
\Delta_s^{\max} = \frac{10}{441.67} \approx 0.0227.
\]

We aim to minimize \( \Delta_s = s_{\max} - s_{\min} \) subject to
\[
\Delta_s \leq \Delta_s^{\max}, \quad s_{\min} \leq s_i \leq s_{\max}, \quad s_i > 0, \quad \text{for } i = 1, 2, 3.
\]

We now solve the optimization problem.

\emph{Case 1: Uniform scaling is feasible.}

If \( \Delta_s^{\max} \geq 0 \), setting \( s_1 = s_2 = s_3 = s \) minimizes \( \Delta_s = 0 \) and satisfies the constraint \( \Delta_s \leq \Delta_s^{\max} \).

\emph{Case 2: Uniform scaling is not feasible.}

If \( \Delta_s = 0 \) does not satisfy the constraint \( d_B(D, D_S) \leq \epsilon \), we set \( \Delta_s = \Delta_s^{\max} \). We choose \( s_{\min} = 1 \) and \( s_{\max} = 1 + \Delta_s^{\max} \).

Assign scaling factors:
- \( s_1 = s_{\min} = 1 \) (e.g., Red channel).
- \( s_2 = s_{\min} = 1 \) (e.g., Green channel).
- \( s_3 = s_{\max} = 1 + \Delta_s^{\max} \) (e.g., Blue channel).

We now verify the topological constraint by using Theorem 1
\[
d_B(D, D_S) \leq \Delta_s \cdot \operatorname{diam}(X) = \Delta_s^{\max} \cdot \operatorname{diam}(X) = \epsilon.
\]
Thus, the topological perturbation remains within the specified tolerance.

We then implement the scaling transformation. We apply the scaling transformation \( S \) to the RGB values of each pixel \( p \)
\[
S(R_p, G_p, B_p) = (s_1 R_p, s_2 G_p, s_3 B_p).
\]

For instance, with \( s_1 = s_2 = 1 \) and \( s_3 = 1 + \Delta_s^{\max} \), the Blue channel is slightly enhanced, introducing variation while preserving the overall color relationships and topology.

The persistence diagrams \( D \) and \( D_S \) capture the topological features of the images before and after scaling, respectively. Features in images often correspond to edges, textures, and regions of uniform color.

By ensuring \( d_B(D, D_S) \leq \epsilon \), we guarantee that the significant topological features (e.g., objects and shapes within the image) are preserved. Minor variations introduced by the scaling are controlled and do not distort the essential structure of the image.

\subsubsection*{Detailed Numerical Example}

Suppose we have an image with the following characteristics:

- Maximum RGB values observed in the image: \( (R_{\max}, G_{\max}, B_{\max}) = (200, 180, 220) \).

- Minimum RGB values observed in the image: \( (R_{\min}, G_{\min}, B_{\min}) = (50, 60, 40) \).

Compute the dataset diameter:
\[
\operatorname{diam}(X) = \sqrt{(200 - 50)^2 + (180 - 60)^2 + (220 - 40)^2} \approx \sqrt{150^2 + 120^2 + 180^2} \approx 263.02.
\]

Given \( \epsilon = 5 \), the maximum allowed scaling variability is
\[
\Delta_s^{\max} = \frac{5}{263.02} \approx 0.019.
\]

Set \( s_{\min} = 1 \) and \( s_{\max} = 1 + 0.019 = 1.019 \).

We assign scaling factors
- \( s_1 = 1 \) (Red channel).

- \( s_2 = 1 \) (Green channel).

- \( s_3 = 1.019 \) (Blue channel).

We now compute the upper bound on \( d_B(D, D_S) \):
\[
d_B(D, D_S) \leq \Delta_s \cdot \operatorname{diam}(X) = 0.019 \times 263.02 \approx 5 \leq \epsilon.
\]

The slight increase in the Blue channel intensifies blue hues in the image without significantly altering the topological features. Edges, contours, and textures remain largely unaffected, ensuring that the augmented image is still suitable for training object recognition models.

\vspace{0.15pt}

\subsection{Case Study: Multimodal Data Normalization}

In many modern machine learning applications, datasets consist of multimodal data, combining features from different sources or modalities, such as text, images, audio, and numerical measurements. These modalities often have inherently different scales and units, which can lead to imbalances in feature importance when training machine learning models. Proper normalization across modalities is crucial to ensure that each feature contributes appropriately to the model's learning process \cite{ghahremani2023regbn}.

Using our mathematical framework, we aim to determine optimal scaling factors for features from each modality to align their scales, minimize scaling variability \( \Delta_s \), and preserve the topological structure of the combined dataset.

\subsubsection*{Objective}

Find scaling factors \( s_i > 0 \) for features across different modalities that minimize \( \Delta_s = s_{\max} - s_{\min} \) while ensuring the topological stability of the multimodal dataset under scaling transformations.

\subsubsection*{Context and Challenges}

Consider a dataset \( X \) comprising features from two modalities:

1. \textbf{Text Features:} Represented using numerical vectors obtained from techniques like word embeddings (e.g., Word2Vec, GloVe) or sentence embeddings. These vectors typically reside in high-dimensional spaces (e.g., \( \mathbb{R}^{300} \)) and have values in a range determined by the embedding method.

2. \textbf{Image Features:} Extracted using convolutional neural networks (CNNs), resulting in feature vectors in \( \mathbb{R}^{n} \), where \( n \) depends on the network architecture and the layer from which features are extracted.

The scales of these features can differ significantly due to the nature of the data and the extraction methods used. If left unnormalized, features from one modality may dominate the learning process, leading to suboptimal model performance.

\subsubsection*{Analysis}

The combined dataset \( X \) consists of feature vectors \( x \in \mathbb{R}^n \), where \( n = n_{\text{text}} + n_{\text{image}} \)
\[
x = (x_{\text{text}}, x_{\text{image}}),
\]
where \( x_{\text{text}} \in \mathbb{R}^{n_{\text{text}}} \) and \( x_{\text{image}} \in \mathbb{R}^{n_{\text{image}}} \).

We then compute the range or variance of each feature to assess the scaling disparity between modalities

- For text features, calculate \( \operatorname{Range}_{\text{text}} = \max_i x_{\text{text}, i} - \min_i x_{\text{text}, i} \).

- For image features, calculate \( \operatorname{Range}_{\text{image}} = \max_i x_{\text{image}, i} - \min_i x_{\text{image}, i} \).

Suppose we find that \( \operatorname{Range}_{\text{text}} \approx 1 \) (e.g., embeddings normalized to unit length), while \( \operatorname{Range}_{\text{image}} \approx 100 \) (e.g., features with larger magnitudes).

We then calculate the diameter of the combined dataset \( X \)
\[
\operatorname{diam}(X) = \max_{p, q \in X} \| p - q \|_2.
\]
Given the disparity in feature scales, the diameter will be dominated by the modality with larger feature ranges (in this case, image features).

We can select a tolerance \( \epsilon > 0 \) representing the maximum acceptable topological perturbation. Compute the maximum allowed scaling variability
\[
\Delta_s^{\max} = \frac{\epsilon}{\operatorname{diam}(X)}.
\]
For instance, if \( \operatorname{diam}(X) = 200 \) and \( \epsilon = 5 \), then:
\[
\Delta_s^{\max} = \frac{5}{200} = 0.025.
\]

We aim to find scaling factors \( s_i > 0 \) for each feature that minimize \( \Delta_s = s_{\max} - s_{\min} \) while ensuring \( \Delta_s \leq \Delta_s^{\max} \). The optimization problem is
\[
\begin{aligned}
& \min_{s_1, \ldots, s_n} && \Delta_s = s_{\max} - s_{\min}, \\
& \text{subject to} && \Delta_s \leq 0.025, \\
& && s_{\min} \leq s_i \leq s_{\max}, \quad \forall i, \\
& && s_i > 0, \quad \forall i.
\end{aligned}
\]

Given the structure of the dataset, we can assign scaling factors based on modality:

- \textbf{Text Features:} Apply a scaling factor \( s_{\text{text}} \) to all text features.

- \textbf{Image Features:} Apply a scaling factor \( s_{\text{image}} \) to all image features.

Our variables reduce to \( s_{\text{text}} \) and \( s_{\text{image}} \), simplifying the problem.

We now solve the optimization problem.

\emph{Case 1: Equalizing the Scales}

Aim to adjust \( s_{\text{text}} \) and \( s_{\text{image}} \) to equalize the ranges of the modalities:

1. Compute the scaling factors required to normalize the ranges:
\[
s_{\text{text}} = \frac{\operatorname{Range}_{\text{image}}}{\operatorname{Range}_{\text{text}}}, \quad s_{\text{image}} = 1.
\]

For example, with \( \operatorname{Range}_{\text{image}} = 100 \) and \( \operatorname{Range}_{\text{text}} = 1 \):
\[
s_{\text{text}} = 100, \quad s_{\text{image}} = 1.
\]

2. Compute \( \Delta_s = s_{\max} - s_{\min} = s_{\text{text}} - s_{\text{image}} = 100 - 1 = 99 \).

3. Check if \( \Delta_s \leq \Delta_s^{\max} \):
\[
99 \leq 0.025 \quad \text{(False)}.
\]
The scaling variability is too large, violating the constraint.

\emph{Case 2: Minimizing \( \Delta_s \) Within Constraints}

Set \( \Delta_s = \Delta_s^{\max} = 0.025 \). Choose \( s_{\min} = 1 \) and \( s_{\max} = 1 + \Delta_s^{\max} = 1.025 \).

Assign scaling factors:
- \( s_{\text{text}} = s_{\min} = 1 \).
- \( s_{\text{image}} = s_{\max} = 1.025 \).

Now, \( \Delta_s = s_{\max} - s_{\min} = 0.025 \), satisfying the constraint.

We now verify the topological constraint by using Theorem 1.
\[
d_B(D, D_S) \leq \Delta_s \cdot \operatorname{diam}(X) = 0.025 \times 200 = 5 \leq \epsilon.
\]
Thus, the topological perturbation remains within the specified tolerance.

We now apply the scaling transformation. Begin scaling the features

- For text features: \( x_{\text{text}}' = s_{\text{text}} \cdot x_{\text{text}} \).

- For image features: \( x_{\text{image}}' = s_{\text{image}} \cdot x_{\text{image}} \).

By adjusting the scaling factors,

- The features from both modalities contribute more equally during model training.

- The topological features of the combined dataset are preserved, preventing distortion of the data's intrinsic structure.

- The model can learn meaningful relationships across modalities without bias toward one modality due to scale differences.

\vspace{0.3pt}

\section{Conclusion}

Throughout the paper, we have shown that the bottleneck distance \( d_B(D, D_S) \) between persistence diagrams under non-uniform scaling \( S \) satisfies:
\[
d_B(D, D_S) \leq \Delta_s \cdot \operatorname{diam}(X),
\]
where \( \Delta_s = s_{\max} - s_{\min} \). This establishes a direct relationship between scaling variability and topological perturbation.

Our results extend to higher homology dimensions \( k \), alternative metrics such as Wasserstein distances \( W_p(D, D_S) \), iterative transformations, and random scaling factors. Specifically, for the \( k \)-th homology, we have:
\[
d_B(D^k, D_S^k) \leq \Delta_s \cdot \operatorname{diam}_k(X),
\]
where \( \operatorname{diam}_k(X) \) is the maximum diameter among \( (k+1) \)-tuples in \( X \).

The proposed framework minimizes \( \Delta_s \) while maintaining \( d_B(D, D_S) \leq \epsilon \), ensuring topological stability. This guarantees that data augmentation via scaling transformations preserves essential features, providing a robust foundation for applications in machine learning and multimodal data analysis.

\vspace{0.3in}

\bibliographystyle{ieeetr}
\bibliography{references}

\end{document}